\documentclass[12pt]{article}

\newcommand{\CC}{{\mathbf C}}
\newcommand{\ZZ}{{\mathbf Z}}
\newcommand{\NN}{{\mathbf N}}
\newcommand{\QQ}{{\mathbf Q}}
\newcommand{\Gal}{\mbox{Gal}}
\newcommand{\mod}{\mbox{mod}}
\newcommand{\Hom}{\mbox{Hom}}

\title{On Kummer extensions of the power series field}
\author{Jos\'e M. Tornero\footnote{Supported by JdA (FQM 218)
and MCyT (BFM2001-3207 and FEDER).}}
 
\date{November, 2002}

\begin{document}

\maketitle

\abstract{In this paper we study the Kummer 
extensions of a power series field $K=k((X_1,...,X_r))$, where $k$ is an 
algebraically closed field of arbitrary characteristic.}

\section{Terminology and notation}

Let $k$ be an algebraically closed field, $X_1,...,X_r$ indeterminates formally 
independent over $k$, and let $K$ and $L_m$ be the fields
$$
K=k\left(\left(X_1,...,X_r\right)\right),
L_m=k\left(\left(X_1^{1/m},...,X_r^{1/m}\right)\right),
$$
where $m$ is a non negative integer, not divisible by the characteristic of $k$.

The extension $K \subset L_m$ is trivially normal, finite and separable, its 
Galois group being $G \simeq (C_m)^r$, where $C_m$ stands for the cyclic group 
of $m$ elements. The elements of $G$ will be noted
\begin{eqnarray*}
\left( a_1,...,a_r \right) : L_m & \longrightarrow & L_m, \quad
0 \leq a_i < m \\
X_l & \longmapsto & \omega^{a_l} X_l
\end{eqnarray*}
where $\omega \in k$ is an $m$--th primitive root of the unity.

Let $R$ and $S_m$ be the rings
$$
R=k\left[\left[X_1,...,X_r\right]\right],\;\;
S_m =k\left[\left[X_1^{1/m},...,X_r^{1/m}\right]\right].
$$
The elements of $S_m$ will be called {\em Puiseux power series}.

Our field of study will be Kummer extensions of $K$. In order to do that, 
recall (\cite{Artin}) that a {\em Kummer extension} of exponent $n$ of a 
field $F$ (which must containing a primitive $n$--th root of the unity and hence 
its characteristic cannot divide $n$), is the splitting field of a polynomial
$$
\left( Z^n - \alpha_1 \right)...\left( Z^n - \alpha_q \right), \mbox{ with }  
\alpha_1,...,\alpha_q \in F.
$$

Our purpose is to prove the following result:

{\bf Theorem.--} {\it
Let $K \subset K'$ be an algebraic separable extension. Then $K'$ can be 
generated by a set of monomials lying in some $S_m$ if and only if there 
exists a Puiseux power series $\zeta \in S_m$ such that $K' = K \left[ 
\zeta \right]$. 
}

It becomes obvious that all separable extensions of $K$ generated by monomials 
lying in some $S_m$ are Kummer extensions; so this proves that all subextensions
of $L_m$ are Kummer.

Notice that an extension generated by monomials in $S_m$ should contain a 
Puiseux power series which generates it (using the primitive element theorem, 
as $k$ must be infinite). In the next section, we will prove the converse.
Finally we will make some general remarks about Kummer
extensions not contained in any $L_m$.

\section{Distinguished exponents of a Puiseux power series}

If $\zeta \in S_m$ is written as
$$
\zeta = \sum c_{i_1...i_r} X_1^{i_1/m}... X_r^{i_r/m}, \; c_{i_1...i_r} 
\in k,
$$
then the set 
$$
\Delta \left( \zeta \right) = \left\{ (i_1,...,i_r) \; | \; 
c_{i_1...i_r} \neq 0 \right\}  \subset \NN^r
$$
will be called (a bit carelessly) the {\em set of exponents} of $\zeta$.

{\bf Definition.--} 
Given $\zeta \in S_m$, a finite subset 
$$
\left\{ \left(i^{(1)}_1,...,i^{(1)}_r\right),...,\left(
i^{(s)}_1,...,i^{(s)}_r \right) \right\} \subset 
\Delta \left( \zeta \right)
$$
will be called a {\em set of distinguished exponents} of $\zeta$ if 
$$
K \left( \zeta \right) = 
K \left( X_1^{i^{(1)}_1/m}...X_r^{i^{(1)}_r/m},...,
X_1^{i^{(s)}_1/m}...X_r^{i^{(s)}_r/m} \right).
$$

In order to prove that all separable extensions of $K$ generated by a Puiseux
power series are Kummer extensions, it suffices to
check that all Puiseux power series in $S_m$, where $m$ is not divisible by
ch$(k)$,
possess a set of distinguished exponents.

We will describe here a process for obtaining such a set for a given
series $\zeta \in S_m$. First of all we fix a total ordering in $\NN^r$, 
say $\prec$, and assume that $m$ is the minimal denominator for $\zeta$
(that is, $\zeta \notin S_q$ for all $q<m$).

For a given matrix $A$ of $t$ rows and $u$ columns, whose elements are 
integers, we will write
$$
(l)\gcd (A) = \gcd \left( \mbox{minors of order $l$ in $A$} \right),
$$
for all $l=1,...,\min\{ t, u \}$.

{\bf Step 1.--} 
Consider the $r \times r$ matrix
$$
M_0 = \left( \begin{array}{cccc} 
m & 0 & ... & 0 \\ 
0 & m & ... & 0 \\
\vdots & \vdots & \ddots & \vdots \\
0 & 0 & ... & m
\end{array} \right),
$$
which obviously verifies $(r)\gcd (M_0) = m^r$.

{\bf Step 2.--} 
Define the sets $\Delta_0 = \Delta (\zeta)$ and
$$
\Delta'_0 = \left\{ \left( i_1,...,i_r \right) \in \Delta_0 \; \left| \;
(r)\gcd (M_0) = (r)\gcd \left( M_0 \; \begin{array}{|c} 
i_1 \\ 
\vdots \\
i_r 
\end{array}
\right) \right. \right\}.
$$
(These exponents are trivially those representing monomials of $\zeta$ which 
lie in $R$).

{\bf Step 3.--} 
Write $\Delta_1 = \Delta_0 \setminus \Delta'_0$, 
define the first distinguished pair by
$$
\left( i^{(1)}_1,...,i^{(1)}_r \right) = \min_\prec \left( \Delta_1 \right);
$$  
and consider the matrix
$$
M_1 = \left( M_0 \; \begin{array}{|c} 
i^{(1)}_1 \\ 
\vdots \\
i^{(1)}_r
\end{array} \right).
$$

{\bf Step 4.--} 
Once the distinguished pairs 
$$
\left( i^{(1)}_1,...,i^{(1)}_r \right),..., 
\left( i^{(l)}_1,...,i^{(l)}_r \right),
$$ 
the set $\Delta_l$ and the matrix $M_l$ are defined, consider
$$
\Delta'_l = \left\{ \left( i_1,...,i_r \right) \in \Delta_l \; \left| \;
(r)\gcd (M_l) = (r)\gcd \left( M_l \; \begin{array}{|c} 
i_1 \\ 
\vdots \\
i_r 
\end{array}
\right) \right. \right\}.
$$

{\bf Step 5.--} 
Write $\Delta_{l+1} = \Delta_l \setminus \Delta'_l$, 
define the $(l+1)$--th distinguished pair by
$$
\left( i^{(l+1)}_1,...,i^{(l+1)}_r \right) = \min_\prec \left( 
\Delta_{l+1} \right);
$$  
and consider the matrix
$$
M_{l+1} = \left( M_l \; \begin{array}{|c} 
i^{(l+1)}_1 \\ 
\vdots \\
i^{(l+1)}_r 
\end{array} 
\right).
$$

{\bf Remark.--} The previous procedure must give a finite number of 
distinguished pairs, as for every $l>0$ we have
$$
(r)\gcd \left( M_{l-1} \right) > (r) \gcd \left( M_l \right),
$$
so we must end up with a finite set 
$$
P = \left\{ \left( i^{(1)}_1,...,i^{(1)}_r \right),...,
\left(i^{(s)}_1,..., i^{(s)}_r \right)\right\}.
$$

From now on we will write for short
$$
K[P] = K \left[ X_1^{i^{(l)}_i}...X_r^{i^{(l)}_r} \; | \; l=1,...,s \right].
$$

Now $K[P] \subset K[\zeta]$, as every element of $G$ leaving $\zeta$ fixed, 
does so with the monomials having exponents in $P$. So, for proving that $P$
is a set of distinguished monomials, it suffices proving the following
result:

{\bf Proposition.--} {\it Let there be
$$
P_1 = \left\{ X_1^{j^{(1)}_1/m}... X_r^{j^{(1)}_r/m},..., X_1^{j^{(t)}_1/m}
... X_r^{j^{(t)}_r/m} \right\},
$$
$$
P_2 = P_1 \cup \left\{  X_1^{j^{(t+1)}_1/m} ...  X_r^{j^{(t+1)}_r/m} \right\}
$$
two sets of monomials in $S_m$ (not in any $S_q$, with $q<m$), such that
$$
(r) \gcd \left( M_1 \right) = (r) \gcd \left( M_2 \right),
$$
where
$$
M_1 = \left( 
\begin{array}{cccccc}
m & ... & 0 & j^{(1)}_1 & ... &  j^{(t)}_1 \\
\vdots & \ddots & \vdots & \vdots & & \vdots \\
0 & ... & m &  j^{(1)}_r & ... &  j^{(t)}_r
\end{array} \right),
$$
$$ 
M_2 = \left( 
\begin{array}{ccccccc}
m & ... & 0 & j^{(1)}_1 & ... &  j^{(t)}_1 &  j^{(t+1)}_1 \\
\vdots & \ddots & \vdots & \vdots & & \vdots & \vdots \\
0 & ... & m &  j^{(1)}_r & ... &  j^{(t)}_r &  j^{(t+1)}_r
\end{array} \right).
$$

Then $K \left[ P_1 \right] = K \left[ P_2 \right]$.}

{\bf Proof.--} The point is proving $K \left[ P_1 \right] \supset K \left[ 
P_2 \right]$ and, for this, it is necessary and sufficient showing that, if
we call $G_k = \Gal \left( L_m / K\left[ P_l \right] \right)$, for $l=1,2$;
then $G_1 = G_2$. 

Define the set
$$
H_1 = \left\{ \left( i_1,...,i_r \right) \in (\ZZ / \ZZ m)^r  \left| 
X_1^{i_1/m} ... X_r^{i_r/m} \in K \left[ P_1 \right] \right. \right\}.
$$

So, $H_1$ contains, up to multiples of $m$ in all coordinates, those 
monomials which remain fixed by the elements of $G_1$. Writing up these 
elements in the form $\left( a_1,...,a_r \right)$ it means that
$$
\left( i_1,...,i_r \right) \in H_1 \; \Longleftrightarrow \;
\sum_{l=1}^r a_l i_l = 0 \; \; (\mod \; m), \; \forall 
\left( a_1,...,a_r \right) \in G_1,
$$
and also, in particular,
$$
H_1 = \left\langle \left( j^{(1)}_1,...,j^{(1)}_r \right),...,
\left( j^{(t)}_1,...,j^{(t)}_r \right) \right\rangle.
$$

Therefore $H_1$ is clearly a subgroup of $G$ (non--canonically identified 
with
$(\ZZ / \ZZ m)^r$), but it also admits another interpretation. In fact,
$$
H_1 \simeq \Hom \left( G/G_1, \ZZ / \ZZ m \right),
$$
identifying $\left( i_1,...,i_r \right) \in H_1$ with
\begin{eqnarray*}
f_{\left( i_1,...,i_r \right)}: G/G_1 & \longrightarrow & \ZZ / \ZZ m \\
\left( x_1,...,x_r \right) + G_1 & \longmapsto & \sum_{l=1}^r x_l i_l 
\end{eqnarray*}

As $G$ is the direct sum of $r$ cyclic groups of order $m$, we have that
$G/G_1$ can be written up as
$$
G/G_1 = C_{a_1} \oplus ... \oplus C_{a_c}, \mbox{ where } a_l | m, \;
\forall l=1,...,c.
$$

This leads to
$$
H_1 \simeq \Hom \left( G/G_1, \ZZ / \ZZ m \right) \simeq 
\bigoplus_{l=1}^c \Hom \left( C_{a_l}, \ZZ / \ZZ m \right) \simeq
\bigoplus_{l=1}^c C_{a_l} \simeq G/G_1,
$$
as $a_l|m$, for all $l$.

On the other hand $\left| G/G_1 \right|$ (that is, $\left[ K\left[ 
P_1 \right] : K \right]$), is precisely $\left| H_1 \right|$ and hence, 
$$
\left| G_1 \right| = \left| G/H_1 \right|.
$$

Let us calculate $\left| G/H_1 \right|$. First of all, instead of writing 
the group as
$$
( \ZZ / \ZZ m )^r / \left\langle \left( j^{(1)}_1,...,j^{(1)}_r \right),...,
\left( j^{(t)}_1,...,j^{(t)}_r \right) \right\rangle,
$$
we will do it as $\ZZ^r / \widehat{H_1}$, where
$$
\widehat{H_1} =  \left\langle (m,0,...,0),...,(0,0,...,m), 
\left( j^{(1)}_1,...,j^{(1)}_r \right),..., \left( j^{(t)}_1,...,j^{(t)}_r 
\right) \right\rangle.
$$

Let us write $\varphi$ a generic element of $\Hom ( \ZZ^r, \QQ / \ZZ)$ with
$\widehat{H_1} \subset \ker(\varphi)$, and $\tilde{\varphi}$ its factorization
through $\ZZ^r / \widehat{H_1}$. According to \cite{Bourbaki}, prop. 8; 
 
$$
\ZZ^r / \widehat{H_1} \simeq \Hom \left( \ZZ^r / \widehat{H_1}, \QQ / \ZZ 
\right),
$$ 
and each of these morphisms is characterized by the images of the 
canonical generating set of $\ZZ^r / \widehat{H_1}$, say
$$
\alpha_l = \tilde{\varphi} \left( e_l + \widehat{H_1} \right),
$$
where $e_l$ stands for the $l$--th element of the canonical basis of 
$\ZZ^r$.

But $\widehat{H_1} \subset \ker \left( \varphi \right)$ is equivalent to
$$
\left( \alpha_1 \;\; ... \;\; \alpha_r \right)
\left( \begin{array}{cccccc}
m & ... & 0 & j^{(1)}_1 & ... &  j^{(t)}_1 \\
\vdots & \ddots & \vdots & \vdots & & \vdots \\
0 & ... & m &  j^{(1)}_r & ... &  j^{(t)}_r
\end{array} \right) = ( 0 \;  ... \; 0).
$$

Again by \cite{Bourbaki}, cor. 1, we can find some linear forms $L_1,...,L_r$
with
 coefficients on $\ZZ$ such that the previous relations are equivalent to
$$
\left( L_1 \left( \alpha_1,...\alpha_r \right) \;\; \ldots \;\; 
L_r \left( \alpha_1 , ... \alpha_r \right) \right) 
\left( \begin{array}{cccccc}
\eta_1 & ... & 0 & 0 & ... & 0 \\
\vdots & \ddots & \vdots & \vdots & & \vdots \\
0 & ... & \eta_r & 0 & ... & 0
\end{array} \right)
= ( 0 \;  ... \; 0),
$$
where $\eta_1 = (1)\gcd \left( M_1 \right)=1$ and $\eta_1 ... \eta_l = 
(l) \gcd \left( M_1 \right)$. Therefore, as 
this equality must hold in $\QQ / \ZZ$, it is plain that there are exactly
$(r) \gcd \left( M_1 \right)$ different morphisms in $\Hom \left( \ZZ^r / 
\widehat{H_1}, \QQ / \ZZ  \right)$, hence
$$
\left| G_1 \right| = \frac{m^r}{\left| H_1 \right|} = \frac{m^2}{m^2/ (r) 
\gcd \left( M_1 \right)} = (r) \gcd \left( M_1 \right).
$$

Doing exactly the same with $G_2$ we find
$$
\left| G_2 \right| = (r) \gcd \left( M_2 \right) = (r) \gcd \left( M_1 
\right) = \left| G_1 \right|.
$$

This finishes the proof, as $G_1 \subset G_2$.

{\bf Corollary.--} {\it If $\zeta \in S_m$, having a set of distinguished
exponents
$$
P= \left\{ \left(i^{(1)}_1,...,i^{(1)}_r\right),...,\left(
i^{(s)}_1,...,i^{(s)}_r \right) \right\} \subset 
\Delta \left( \zeta \right),
$$
then
$$
\Delta \left( \zeta \right) \subset \sum_{l=i}^s \ZZ \left(i^{(l)}_1,
...,i^{(l)}_r \right), \;\; \mbox{\rm mod } \ZZ m \times ... \times \ZZ m
$$
}

{\bf Remark.--} This process enables to compute (on equal footing) 
two well--known (sets of) 
arithmetic data which are most useful in algebraic geometry, as are the 
Puiseux pairs of a plane curve (\cite{Zariski}) and the characteristic 
pairs of a quasi--ordinary surface (\cite{Lipman}) (both of them for $k=\CC$).

We will do the Puiseux pairs case. So, assume we have a plane algebroid 
curve given by $f(X,Y) \in  \CC [[X,Y]]$ and a Puiseux branch, which can 
always be represented (up to a change of variables) as
$$
\begin{array}{rcl}
Y = \zeta \left( X^{1/m} \right) & = & \displaystyle c_{\beta_1} X^{\beta_1/m} + 
\sum_{l=1}^{h_1} c_{\beta_1+le_1} X^{\left( \beta_1 +le_1 \right) / m} 
+ \; ...   \\ \\
& & \displaystyle \quad ... \; + c_{\beta_g} X^{\beta_g/m} + \sum_{l=1}^\infty
c_{\beta_g+le_g} X^{\left( \beta_g +le_g \right) / m},
\end{array}
$$
where we can assume $m<\beta_1<...<\beta_g$, $\beta_k \notin \ZZ m$ for
all $k=1,...,g$ and, in addition, if we call
$$
\beta_1 = p_1 e_1, \; m = q_1e_1, \; \gcd \left( p_1, q_1 \right) = 1
$$
$$
e_{l-1} = q_le_l, \; \beta_l = p_le_l, \;
\gcd \left( p_l, q_l \right) = 1; \; \; \forall \, l=2,...,g,
$$
then the pairs $\left( p_1,q_1 \right),...,\left( p_g,q_g \right)$ are 
called the {\em Puiseux pairs} of the curve. Note that these pairs are
determined (and they determine as well) by the set
$$
\left\{ m,\beta_1,...,\beta_g \right\},
$$
called by Zariski the {\em characteristic of the branch $\zeta$}. Also is 
direct from the formulae above that
$$
e_1 = \gcd \left( m, \beta_1 \right), \; \; e_l = \gcd \left( e_{l-1},
\beta_l \right), \; \; \forall \, l=2,...,g.
$$

If we apply our process to the set of exponents on $\Delta \left( \zeta
\right)$ using, for instance, the natural ordering on $\NN$, we start up 
with
$$
M_0 = ( m )
$$
and then choose the smaller element on $\Delta$, that is, $\beta_1$, which,
by the above conditions, happens to verify $ \gcd \left( m, \beta_1 \right) 
< m$, so $i^{(1)} = \beta_1$.

Assume we have already computed the first $l$ distinguished exponents, 
which coincide with $\beta_1,...,\beta_l$ (necessarily in this order because
of our choosing of the ordering on $\ZZ$). Then we have the matrix
$$
M_l = \left( m \; \; \beta_1 \; \; ... \; \; \beta_l \right),
$$
and $\gcd \left( m, \beta_1, ... , \beta_l \right) = e_l$, by the above
considerations. We have discarded in previous steps those monomials which 
can be written as a combination of some $\beta_t$ and $e_t$, for $t<l$. In the 
same way, then, we discard now those elements in $\Delta$ which do not
make smaller the previous  $\gcd$, which are, precisely, those which can be
written up as a combination of $\beta_l$ and $e_l$.

By definition of $\beta_{l+1}$, it has to be the minimal element not yet
discarded, and this proves that our procedure must end up computing the
set $\left\{ \beta_1,...,\beta_g \right\}$.

The quasi--ordinary surface case is similar; in fact there are no 
substantial differences between the two cases. However, note that in order to
obtain the characteristic monomials, one has to
take into account that the chosen total ordering must be graded,
as the characteristic monomials of a quasi--ordinary branch $\zeta$ (up to
normalization) are determined by the following facts (\cite{Lipman}):

(1) They are the minimal elements of $\Delta (\zeta)$ for the natural partial
order.

(2) They generate (irredundantly) the same extension field than the branch
itself. 

\section{General Kummer extensions}

We will make now some remarks about Kummer extensions of $K$ of any kind. In order 
to do that observe that we can reduce the problem to that of the splitting 
field of a polynomial $F(Z) = Z^n - \zeta$, where $\zeta \in R$. 

{\bf Remark.--} First of all, mind that, if $\zeta$ is an irreducible power 
series not associated with any of the $X_i$, we cannot hope the splitting field 
of $F$ to be a subextension of some $L_m$. In fact, $\zeta$ defines a valuation
$v_\zeta$ of $K$ that is unramified over $L_m$, as $v_\zeta \left( X_i \right)=0$,
for all $i$,
but it is obviously ramified over the splitting field. However, we can prove 
a resembling result.

{\bf Proposition.--} In the above situation, there are $\alpha_1,...,\alpha_r 
\in \NN$ such that the splitting field of $F$ is a subextension of
$$
K \left( \left( X_1^{\alpha_1/n} \right) \right)...\left( \left( X_r^{\alpha_r/n} 
\right) \right).
$$

{\bf Proof.--} We will do the proof by induction on $r$, being the case $r=1$ 
direct from the so--called Newton--Puiseux Theorem. So assume that, for all 
$\eta \in R$, there exists a set of positive integers 
$\left\{ \alpha_1,...,\alpha_r \right\}$ such that 
$$
K \left[ \sqrt[n]{\eta} \right] = K \left( \left( X_1^{\alpha_1/n} \right)\right)
... \left( \left( X_r^{\alpha_r/n} \right) \right);
$$
and fix a power series $\zeta \in R \left[\left[ X_{r+1} \right] \right]$, with
$\nu (\zeta) = \lambda_0 \geq 0$, where $\nu$ is the usual order with respect to 
$X_{r+1}$. We want to find a root of $Z^n - \zeta \in K'[Z]$, where $K' = k \left( 
\left( X_1,...,X_{r+1} \right)\right)$.

Let us call from now on $\zeta_i$ the approximate $n$--th root of $\zeta$ (up to
order $i$ in $X_{r+1}$), which will be constructed in what follows. The term with 
minimal degree on $X_{r+1}$ of a $\sqrt[n]{\zeta}$, must be of the form $c_{\lambda_0/n} 
X_{r+1}^{\lambda_0/n}$, where it must hold 
$$
c_{\lambda_0/n}^n = a_{\lambda_0},
$$

So there must be a set of monomials $\left\{ X_1^{\alpha_1/n},...,X_r^{\alpha_r/n} 
\right\} \subset S_n$ such that  
$$
K' \left[ \sqrt[n]{a_{\lambda_0}} \right] = K' \left[ c_{\lambda_0/n} \right] \subset 
K' \left( \left( X_1^{\alpha_1/n} \right) \right)... \left( \left( X_r^{\alpha_r/n} 
\right) \right). 
$$

Hence we can write $\zeta_{\lambda_0} = c_{\lambda_0/n} 
X_{r+1}^{\lambda_0/n}$, which verifies
$$
\zeta_{\lambda_0} \in K' \left( \left( X_1^{\alpha_1/n} \right) \right)
... \left( \left( X_r^{\alpha_r/n} \right) \right) \left[ \left[ X_{r+1}^{\lambda_0/n}
\right] \right], \; \; \nu \left( \zeta_{\lambda_0}^n - \zeta \right) = \lambda_1 > 
\lambda_0.
$$

Now, in the same way, the following term with minimal degree on $X_{r+1}$ of
any $\sqrt[n]{\zeta}$, must be of the form $c_{\lambda_1-[\lambda_0(n-1)/n]}
X_{r+1}^{\lambda_1-[\lambda_0(n-1)/n]}$, where it must hold now
$$ 
n c_{\lambda_1-[\lambda_0(n-1)/n]} c_{\lambda_0/n}^{(n-1)/n} =
a_{\lambda_1},  
$$
with $a_{\lambda_1}$ the initial form (with respect to $X_{r+1}$) of
$\zeta_{\lambda_0}^n - \zeta$ and hence
$$
c_{\lambda_1-[\lambda_0(n-1)/n]} \in K' \left( \left( X_1^{\alpha_1/n} 
\right) \right)... \left( \left( X_r^{\alpha_r/n} \right) \right).
$$

We write $\zeta_{\lambda_1} = c_{\lambda_0/n} X_{r+1}^{\lambda_0/n} +
c_{\lambda_1-[\lambda_0(n-1)/n]} X_{r+1}^{\lambda_1-[\lambda_0(n-1)/n]}$
and, as above, it holds
$$
\zeta_{\lambda_1} \in K' \left( \left( X_1^{\alpha_1/n} \right) \right)
... \left( \left( X_r^{\alpha_r/n} \right) \right) \left[ \left[ X_{r+1}^{\lambda_0/n}
\right] \right], \; \; \nu \left( \zeta_{\lambda_1}^n - \zeta \right) = \lambda_2 > 
\lambda_1.
$$
 
Note that, though $\zeta_{\lambda_1}$ is a Puiseux power series, all the 
exponents in $X_{r+1}$ of $\zeta_{\lambda_1}^n$ are positive integers: in fact, 
they are $\{ \lambda_0,\lambda_1,2\lambda_1 - \lambda_0,...,n\lambda_1- (n-1)
\lambda_0 \}$. In particular, this shows $\lambda_2 \in \NN$.

Assume now we have constructed $\zeta_{\lambda_s}$, verifying

\begin{itemize}
\item $\zeta_{\lambda_s} \in K' \left( \left( X_1^{\alpha_1/n} \right) \right)
... \left( \left( X_r^{\alpha_r/n} \right) \right) \left[ \left[ 
X_{r+1}^{\lambda_0/n} \right] \right]$,

\item $\nu \left( \zeta_{\lambda_s}^n - \zeta \right) = \lambda_{s+1} > 
\lambda_s > ... > \lambda_0$, with all $\lambda_i \in \NN$.

\item The series $\zeta_{\lambda_s}$ has the form
$$
\zeta_{\lambda_s} = \sum_{j=0}^s c_{L_j \left( \lambda_1,...,\lambda_j \right) 
- \alpha_j \lambda_0 + \lambda_0/n} X_{r+1}^{L_j \left( \lambda_1,...,\lambda_j 
\right) - \alpha_j \lambda_0 + \lambda_0/n},
$$
where $L_j$ is a linear form with possitive coefficients and $\alpha_j$ a positive 
integer, for $j=1,...,s$; being $L_0=\alpha_0=0$.
\end{itemize}

If we call $a_{\lambda_{s+1}}$ the initial form of $\zeta_{\lambda_s}^n - \zeta$
then it is clear that we have to define
$$
\zeta_{\lambda_{s+1}} = \zeta_{\lambda_s} + c_{\lambda_{s+1}-[\lambda_0(n-1)/n]}
X_{r+1}^{\lambda_{s+1}-[\lambda_0(n-1)/n]},
$$
where it must hold
$$ 
n c_{\lambda_{s+1}-[\lambda_0(n-1)/n]} c_{\lambda_0/n}^{(n-1)/n} =
a_{\lambda_{s+1}},  
$$
and therefore
$$
c_{\lambda_{s+1}-[\lambda_0(n-1)/n]} \in K' \left( \left( X_1^{\alpha_1/n} 
\right) \right)... \left( \left( X_r^{\alpha_r/n} \right) \right).
$$

Moreover, it is plain that $\nu \left( \zeta^n_{\lambda_{s+1}} - \zeta \right)
= \lambda_{s+2} > \lambda_{s+1}$. Finally, note that all the exponents in $X_{r+1}$
appearing in the developement of $\zeta_{\lambda_{s+1}}^n$ are of the type
$$
\sum_{j=0}^s i_j \left[ L_j \left( \lambda_1,...,\lambda_j \right) - \alpha_j \lambda_0
+ \frac{\lambda_0}{n} \right] + i_{s+1} \left( \lambda_{s+1} - \lambda_0 + 
\frac{\lambda_0}{n} \right),
$$
with $i_0+...+i_{s+1}=n$. This implies that all these exponents are positive integers
of the form $L \left( \lambda_1,...,\lambda_{s+1} \right) - \beta \lambda_0$,
with $L$ a linear form with possitive coeffients, $\beta \in \NN$.

In this way it is shown that
$$
\sqrt[n]{\zeta} \in K' \left( \left( X_1^{\alpha_1/n} \right) \right)
... \left( \left( X_r^{\alpha_r/n} \right) \right) \left[ \left[ 
X_{r+1}^{\lambda_0/n} \right] \right].
$$

This finishes the proof of the proposition.

\section{Final comments}

All the arguments given here can be completely translated word--by--word to 
the analytic context. We hope that this work will be useful as a step
to understand the geometry of algebroid (analytic) hypersurfaces which
admit a Puiseux--like parametrization. In fact, distinguished exponents
have proved to be a useful tool for the surface case (characteristic 0), as
shown in \cite{T}. This results have led us to expect that some
deeper application of class field theory tools may help to the study of the
geometry and the topology of these varieties.

\vspace{.5cm}

Jos\'e M. Tornero

{\sc Departamento de \'Algebra, Facultad de Matem\'aticas.

Avda. Reina Mercedes, s/n. 41012 Sevilla (SPAIN)}

E--mail: tornero@algebra.us.es


\begin{thebibliography}{99}

\bibitem{Artin}
{\it E. Artin}, Galois theory. Notre Dame Mathematical 
Lectures, {\bf 2}. University of Notre Dame Press, 1959. 

\bibitem{Bourbaki}
{\it N. Bourbaki}, Alg\`ebre. Ch. VII: Modules sur les anneaux principaux. 
Actualit\'es Sci. Ind., {\bf 1179}. Hermann, 1959.

\bibitem{Lipman}
{\it J. Lipman}, Quasi--ordinary singularities of embedded surfaces. 
Ph. D. Thesis, Harvard University.

\bibitem{T}
{\it J.M. Tornero}, Aspectos locales de singularidades de superficies:
Superificies de Puiseux. Ph. D. Thesis, University of Sevilla.

\bibitem{Zariski}
{\it O. Zariski}, Studies in equisingularity I. Amer. J. Math. {\bf 87} (1965),
507--536.

\end{thebibliography}
\end{document}